\documentclass[a4paper,11pt]{article}
\usepackage{amsmath}
\usepackage{amssymb}
\usepackage{amsmath}
\usepackage{mathrsfs}
\usepackage{amsmath,amssymb,amsthm,latexsym,amscd,mathrsfs}
\usepackage{indentfirst}
\usepackage{graphicx}
\usepackage{booktabs}
\usepackage{array}
\usepackage{chemarrow}
\usepackage{hyperref}
\usepackage{graphicx}
\usepackage[perpage,symbol*]{footmisc}
\usepackage{float}

\numberwithin{equation}{section} \numberwithin{theorem}{section}
\numberwithin{lemma}{section} \numberwithin{corollary}{section}
\numberwithin{definition}{section}
\numberwithin{proposition}{section} \numberwithin{remark}{section}
\numberwithin{example}{section}

\makeatletter
\def\@makefnmark{}

\makeatother

\begin{document}
\vspace{20mm}

\noindent
{\it Translator's Note: } This is an unofficial translation of the original paper
which was written in German. All references should be made to the original paper.
\vspace{8mm}

\noindent
Translation of ~{\it Cliffordalgebren und neue isoparametrische Hyperfl\"{a}chen},~
Math. Z. \textbf{177}, 479--502 (1981) by Dirk Ferus, Hermann Karcher and Hans-Friedrich
M\"{u}nzner.
\vspace{5mm}

\noindent
1. Translated by Thomas E. Cecil, Department of Mathematics and Computer Science, College
of the Holy Cross, Worcester, MA 01610, USA; 

E-mail address: \textbf{cecil@mathcs.holycross.edu}
\vspace{1mm}

\noindent
2. Typed by Wenjiao Yan, School of Mathematical Sciences, Laboratory of Mathematics and Complex Systems, Beijing Normal
University, Beijing 100875, China. 

E-mail address: \textbf{wjyan@mail.bnu.edu.cn}\\

\title{Clifford Algebras and New Isoparametric Hypersurfaces}
\author{Dirk Ferus$^1$, Hermann Karcher$^2$, Hans-Friedrich M\"{u}nzner$^3$
}
\date{}\maketitle

\begin{quotation}
\noindent
1. Fachbereich Mathematik der Technischen Universit\"{a}t, Stra{\ss}e des 17. Juni 135,
D-1000 Berlin 12, Bundesrepublik Deutschland\\
2. Mathematisches Institut der Universit\"{a}t, Wegelerstra{\ss}e 10, D-5300 Bonn,
Bundesrepublik Deutschland\\
3. Studienbereich Mathematik der Universit\"{a}t, Kufsteinerstra{\ss}e, D-2800 Bremen 33,
Bundesrepublik Deutschland\\
\end{quotation}

\noindent
\textbf{1. Introduction.}
\vspace{4mm}

\noindent
Isoparametric hypersurfaces, \emph{i.e.} those with constant principal curvatures, are of special interest
because of their comparatively simple geometry. While in spaces of non-positive constant curvature,
they reduce to the examples of tubes over totally geodesic submanifolds, the theory of isoparametric
hypersurfaces in spheres is rich: the homogeneous examples are completely known through the work of
Hsiang and Lawson \cite{HL}, moreover Takeuchi and Ozeki \cite{OT} found two infinite series of
inhomogeneous examples with $g=4$ distinct principal curvatures. In the following work, by means of
representations of Clifford algebras, for each natural number $m_1$, we will construct an infinite
series of isoparametric hypersurfaces with $g=4$ distinct prinicpal curvatures with multiplicities
$(m_1, m_2, m_1, m_2)$ where $m_2$ grows monotonically in the series. Our examples include all previously
known homogeneous and inhomogeneous $(g=4)$ examples with the exception of two individual homogeneous
cases of dimension $8$ and $18$. Particular attention will be given to the case $m_1\equiv 0~ mod~ 4$.
For this case, as the dimension increases there are more and more families with the same multiplicities
$(m_1, m_2, m_1, m_2)$ and more and more non-isometric compact Riemannian manifolds of the same dimension,
which (modulo an isometry of the tangent space) have the same curvature tensor pointwise. The incongruence
of families with the same multiplicities -- as well as their inhomogeneity in many cases -- can be seen through
the second fundamental forms of the focal submanifolds, which are quite accessible in our representation of
the families.

On isoparametric hypersurfaces, the curvature tensor, as a field of endomorphisms on bivectors
also has constant eigenvalues.
Our examples (as well as the series of Takeuchi and Ozeki) provide a large number of Riemannian manifolds with
these properties. This construction by means of Clifford representations carries over to spaces with indefinite
scalar product and provides examples of isoparametric hypersurfaces in Lorentz spaces of constant curvature which
Nomizu has recently begun to study.

Topologically, our examples of multiplicities $(m_1, m_2, m_1, m_2)$ have the form of $m_2$-sphere bundles over
an $(m_1+m_2)$-sphere bundle over an $m_1$-sphere.

In the two following sections, we give a brief summary of the necessary differential geometric and algebraic foundations.
Section $4$ contains the definition of the new examples; and in the rest we study their geometry.
\vspace{4mm}

\noindent
\textbf{2. Isoparametric hypersurfaces in spheres.}
\vspace{4mm}

In this section, we put together the fundamental properties of isoparametric hypersurfaces which are important to us.
For further details see \cite{Mun}.

\vspace{2mm}

\noindent
\textbf{2.1~Definitions:} An oriented hypersurface in a sphere with constant principal curvatures is called an \emph{isoparametric
hypersurface}. ( If one only requires local constancy, then this definition makes sense for non-orientable hypersurfaces as well,
but~--~as can be shown~--~not globally: so the more general definition is furnished only for orientable hypersurfaces. )
\vspace{2mm}

\noindent
\textbf{2.2~Parallel surfaces:} Given an oriented hypersurface $M$ in the sphere $S^{n+1}$ with principal curvatures $\cot(\phi_i),~0<\phi_i<\pi$
with multiplicities $m_i$, then the parallel hypersurface $M_\epsilon$ at a distance $\epsilon$ in the direction of the normal
has principal curvatures $\cot(\phi_i-\epsilon)$ with the same multiplicities, and $M$ and $M_{\epsilon}$ have the same normal great circles.
The parallel hypersurfaces obtained this way are likewise isoparametric; they form an \emph{isoparametric family.}

The mean curvature $\frac{1}{n}\sum m_i\cot(\phi_i-\epsilon)$ as an analytic function of $\epsilon$ determines the $\phi_i$
through its poles (modulo $\pi$). If all the hypersurfaces in a parallel family have constant mean curvature, then they have constant
principal curvatures and so are isoparametric.
\vspace{2mm}

\noindent
\textbf{2.3~Principal foliations ( curvature foliations. )} The $m_i$-dimensional eigendistribution of the second
fundamental form ( Weingarten map ) $S$ of an isoparametric hypersurface corresponding to the principal curvature $\cot(\phi_i)$ is
integrable and autoparallel. It integrates to a totally geodesic foliation of $M$ by $m_i$-spheres of radius $\phi_i$ in $S^{n+1}$.
The parallel surface at distance $\phi_i$ is a {\it focal manifold} of the isoparametric family, and the hypersurfaces of the
family are tubes over the focal manifolds.
\vspace{2mm}

\noindent
\textbf{2.4~Distribution of values of the principal curvatures.} Consider the principal curvature $\cot \phi_i$ of some
hypersurface of the isoparametric family. Its curvature leaf is a sphere and its center is a point of the ``corresponding"
focal submanifold, with the radii of the sphere as normals. The eigenvalues of the $2^d$ fundamental tensors are
$\cot(\phi_j-\phi_i)$, $j\neq i$, independent of the point and the
normal direction. From this, one gets through an argument not given here, that the principal curvatures of the hypersurfaces
are of the form $\cot(\alpha+\frac{g-i}{g}\pi)$ with multiplicities $m_i$. Here, $g$ is the number of distinct principal curvatures,
$1\leq i\leq g$, and $\alpha\in (0, \pi)$. Moreover, it turns out that $m_i=m_{i+2}$ with indices modulo $g$, see M\"{u}nzner \cite{Mun}.
\vspace{2mm}

\noindent
\textbf{2.5~Isoparametric functions: } If $f:S^{n+1}\rightarrow \mathbb{R}$ is a function whose gradient on each level
surface has constant length, then the regular levels are parallel hypersurfaces with $\texttt{grad}~ f \Big/\|\texttt{grad}~ f\|$ as unit normal field.
The corresponding second fundamental tensor is given by
\begin{equation*}
SV=-\nabla_V\frac{\texttt{grad}~ f}{\|\texttt{grad}~ f\|}=\frac{-1}{\|\texttt{grad}~f\|}\nabla_V\texttt{grad}~ f ;
\end{equation*}
the mean curvature is $-\frac{1}{n}\frac{\Delta f}{\|\texttt{grad}~f\|}$.

Thus the levels form an isoparametric family precisely when $\Delta f$ is constant on the levels, that is, when both of ``the first two
differential parameter" functions are functions of $f$ itself:
$$\|\texttt{grad}~f\|^2=a(f),~~\Delta f =b(f).$$
Whence the name ``isoparametric".
\vspace{2mm}

\noindent
\textbf{2.6~Differential equations of Cartan and M\"{u}nzner:} Conversely, given an isoparametric family and a function $f$
of the oriented distance from a fixed hypersurface of the family, then $f$ obviously satisfies the above differential equations.
The information in $2.4$ can be exploited to specify the right hand side of the differential equations further. Through Cartan and
M\"{u}nzner
$$f=\cot gt,$$
where $t$ is the spherical distance to a focal manifold and as usual $g$ is the number of distinct principal
curvatures, one is led to a particularly good standardization: extend $f$ to a homogeneous function $F$ of degree $g$ on a cone in
$\mathbb{R}^{n+2}$, then the differential equations for $F$ ( with Euclidean differential operators ) are written as:
$$\|\texttt{grad}~F\|^2=g^2r^{2g-2},~~\Delta~F=cr^{g-2},$$ with $c=g^2\frac{m_2-m_1}{2}$ ( $=0$ for odd $g$ ) and $r(x)=\|x\|$.

One sees as a simple consequence that $F$ is a homogeneous polynomial of degree $g$, see \cite{Mun}. In particular, one can extend
each piece of isoparametric hypersurface to a compact algebraic hypersurface.
\vspace{2mm}

\noindent
\textbf{2.7~Rigidity.} Isoparametric hypersurfaces are rigid: each isometry of such a hypersurface can be extended to an
isometry of the sphere; hence the parallel hypersurfaces are mapped to themselves, thus causing an isometry of the whole family.
For hypersurfaces with $g\geq 4$, which are of interest to us here, this follows from the classical rigidity theorem since the
rank of the second fundamental tensor is obviously $\geq 3$. For $g\leq 3$, one can make use of the classification of Cartan \cite{Car}
or otherwise carry out a direct proof.
\vspace{2mm}

\noindent
\textbf{2.8~Minimality.} From the facts collected in $2.4$ one obtains that the focal manifolds of an isoparametric family
are minimal submanifolds of the sphere and each such family contains precisely one minimal hypersurface.
\vspace{4mm}

\noindent
\textbf{3. Clifford systems.}
\vspace{4mm}

Our examples of isoparametric hypersurfaces are constructed from from so-called Clifford systems whose definitions and properties
are collected in this section.
\vspace{2mm}

\noindent
\textbf{3.1~Notation.} For a Euclidean vector space $V$, let $\mathfrak{h}(V)$ \emph{resp.} $\mathcal{O}(V)$ \emph{resp.} $O(V)$ denote the symmetric \emph{resp.} skew-symmetric \emph{resp.} orthogonal endomorphisms of $V$. We let $\langle~,~\rangle$ denote the canonical scalar product on $\mathbb{R}^n$, and on $\mathfrak{h}(V)$
$$\langle A, B\rangle:=\frac{1}{\dim V}\texttt{Trace}(AB).$$
$E_{\pm}(A)$ denotes the eigenspace of the eigenvalue $\pm1$.
\vspace{2mm}

\noindent
\textbf{3.2~Definition.} Let $l,m$ be positive natural numbers.
\begin{enumerate}
  \item The $(m+1)$-tuple $(P_0,...,P_m)$ with $P_i\in \mathfrak{h}(\mathbb{R}^{2l})$ is called a \emph{(symmetric) Clifford system} on $\mathbb{R}^{2l}$  if for all $i, j \in \{0,...,m\}$ we have
      $$P_iP_j+P_jP_i=2\delta_{ij}Id.$$
  \item Let $(P_0,...,P_m)$ and $(Q_0,...,Q_m)$ be Clifford systems on $\mathbb{R}^{2l}$ resp. $\mathbb{R}^{2n}$, then
  $(P_0\oplus Q_0,...,P_m\oplus Q_m)$ is a Clifford system on $\mathbb{R}^{2(l+n)}$, the so-called \emph{direct sum} of
  $(P_0,...,P_m)$ and $(Q_0,...,Q_m)$.
  \item A Clifford system $(P_0,...,P_m)$ on $\mathbb{R}^{2l}$ is called \emph{irreducible} when it is not possible to write $\mathbb{R}^{2l}$ as a direct sum of two positive dimensional subspaces which are invariant under all of the $P_i$.
\end{enumerate}
\vspace{2mm}

\noindent
\textbf{3.3~Representations of Clifford algebras and Clifford systems.} Each representation of a Clifford algebra $\mathcal{C}_{m-1}$
on $\mathbb{R}^l$ can be characterized by $E_1,...,E_{m-1}\in \mathcal{O}(\mathbb{R}^l)$ with
$$E_iE_j+E_jE_i=-2\delta_{ij}Id.$$
One now defines $(P_0,...,P_m)\in \mathfrak{h}(\mathbb{R}^{2l})$ by
$$P_0(u,v):=(u,-v),~P_1(u,v):=(v,u),~P_{1+i}(u,v):=(E_iv,-E_iu),~u,v\in \mathbb{R}^l,$$
then $(P_0,...,P_m)$ is a symmetric Clifford system on $\mathbb{R}^{2l}$ that is irreducible precisely when the representation of
$\mathcal{C}_{m-1}$ is irreducible.

Conversely, one obtains all symmetric Clifford systems in this way: let $(P_0,...,P_m)$ be a Clifford system on $\mathbb{R}^{2l}$,
then since $P_0^2=Id$, $P_0$ has eigenvalues $\pm 1$, and since $P_0P_1+P_1P_0=0$, $P_1$ interchanges the corresponding eigenspaces
$E_{\pm}(P_0)$. Thus $E_+(P_0)$ is of dimension $l$.

Further, $E_+(P_0)$ is invariant under $P_0P_{i+1}, i\in \{1,...,m-1\}$. One can identify $E_+(P_0)$ with $\mathbb{R}^l$, and
thus define through
$$E_i:=P_1P_{i+1}|_{E_+(P_0)}$$
a representation of $\mathcal{C}_{m-1}$ on $\mathbb{R}^l$.
\vspace{2mm}

\noindent
\textbf{3.4~Definition.} Two Clifford systems $(P_0,...,P_m)$ and $(Q_0,...,Q_m)$ on $\mathbb{R}^{2l}$ are called
\emph{algebraically equivalent} if there exists $A\in O(\mathbb{R}^{2l})$ such that $Q_i=AP_iA^t$ for all $i\in\{0,...,m\}$. They
are called \emph{geometrically equivalent} when there exists $B\in O(Span\{P_0,...,$ $P_m\})\subset \mathfrak{h}(\mathbb{R}^{2l})$ such that
$(Q_0,...,Q_m)$ and $(B(P_0),...,B(P_m))$ are algebraically equivalent.
\vspace{2mm}

\noindent
\textbf{3.5~Representation theory.}  With the help of $3.3$, we obtain from the representation theory of Clifford algebras,
see \cite{Hus}, the following results: each Clifford system is algebraically equivalent to a direct sum of irreducible Clifford systems.
An \emph{irreducible} Clifford system $(P_0,...,P_m)$ on $\mathbb{R}^{2l}$ exists precisely for the following values of $m$ and $l=\delta(m)$:

\begin{center}
\begin{tabular}{|c|c|c|c|c|c|c|c|c|c|}
\hline
$m$ & 1 & 2 & 3 & 4 & 5 & 6 & 7 & 8 & $\cdots$ $m$+8 \\
\hline
$\delta(m)$ & 1 & 2 & 4 & 4 & 8 & 8 & 8 & 8 & $\cdots$ 16$\delta(m)$\\
\hline
\end{tabular}
\end{center}

For $m\not \equiv 0 (mod~4)$, there exists exactly one algebraic ( and therefore exactly one geometric ) equivalence class of irreducible
systems. Thus for each positive integer $k$ there exists exactly one algebraic ( or geometric ) equivalence class of Clifford systems
$(P_0,...,P_m)$ on $\mathbb{R}^{2l}$ with $l=k\delta(m)$.

For $m\equiv 0 (mod~4)$, there exist exactly two algebraic classes of irreducible systems. These are distinguished from each other by
the choice of sign in
$$\texttt{Trace}(P_0\cdots P_m)=\pm \texttt{Trace} ~Id =\pm 2 \delta(m).$$
Hence (by replacing $P_0$ by $-P_0$) there exists exactly one geometric equivalence class in this case also.

The absolute trace $$|\texttt{Trace}(P_0\cdots P_m)|$$ is obviously an invariant under geometric equivalence. If one constructs
all direct sums of both of the irreducible algebraic classes with altogether $k$ summands, then this invariant takes on $[\frac{k}{2}]+1$
different values. Thus for $m\equiv 0 (mod~4)$, there are exactly $[\frac{k}{2}]+1$ geometric equivalence class of Clifford systems
on $\mathbb{R}^{2l}$ with $l=k\delta(m)$.
\vspace{2mm}

\noindent
\textbf{3.6~Definition.} Let $(P_0,...,P_m)$ be a Clifford system on $\mathbb{R}^{2l}$. The unit sphere in $Span\{P_0,...,P_m\}\subset \mathfrak{h}(\mathbb{R}^{2l})$ is called the \emph{Clifford sphere} determined by the system and is denoted
$\Sigma(P_0,...,P_m)$.
\vspace{2mm}

\noindent
\textbf{3.7~Properties of the Clifford sphere.} Our construction will not depend on $(P_0,...,P_m)$ but only on
$\Sigma(P_0,...,P_m)$. We therefore begin with several important properties of Clifford spheres.

$(i)$ For each $P\in \Sigma(P_0,...,P_m)$, we have $P^2=Id$. Conversely, if $\Sigma$ is a unit sphere in a linear subspace $\mathbb{R}\Sigma\subset \mathfrak{h}(\mathbb{R}^{2l})$ such that $P^2=Id$ for all $P\in \Sigma$, then every orthonormal basis of $\mathbb{R}\Sigma$ is a Clifford system on $\mathbb{R}^{2l}.$

$(ii)$ Two Clifford systems on $\mathbb{R}^{2l}$ are geometrically equivalent if and only if their Clifford spheres are conjugate to one
another under an orthogonal transformation of $\mathbb{R}^{2l}.$

$(iii)$ The function
$$H(x)=\sum_{i=0}^{m}\langle P_ix,x\rangle^2$$
depends only on $\Sigma(P_0,...,P_m)$, and not on the choice of orthonormal basis $(P_0,...,P_m)$. For $P\in \Sigma(P_0,...,P_m)$, we have
$H(Px)=H(x)$ for all $x$. (To prove this, choose $P_0,...,P_m$ orthonormal with $P_0=P$.)

$(iv)$ For orthonormal $Q_1,..,Q_r\in \Sigma(P_0,...,P_m)$ since $Q_iQ_j=-Q_jQ_i$ for all $i\neq j$, we have
\begin{eqnarray*}
&&Q_1\cdots Q_r\in \mathfrak{h}(\mathbb{R}^{2l}),\quad for~ r\equiv 0,1 (mod~4)\\
&&Q_1\cdots Q_r\in \mathcal{O}(\mathbb{R}^{2l}),\quad for~ r\equiv 2,3 (mod~4).
\end{eqnarray*}

Further, $Q_1\cdots Q_r$ is uniquely determined by an orientation of $Span(Q_1,...,Q_r)$. $SO(r)$ is generated as a group by rotations
of two-dimensional coordinate planes, and since one can bring any two $Q_i$'s next to each other through permutation modulo signs,
it suffices to do the proof for $r=2$. This is an easy direct calculation.

$(v)$ For $P,Q\in Span\{P_0,...,P_m\}$ and $x\in \mathbb{R}^{2l}$, we have
$$\langle Px,Qx\rangle=\langle P,Q\rangle\langle x,x\rangle.$$
\vspace{4mm}

\noindent
\textbf{4. The new examples.} We now give the new series with $g=4$ distinct principal curvatures and study their geometry.
\vspace{2mm}

\noindent
\textbf{4.1}~
\noindent
\textbf{Theorem:}\,\,
{\itshape
Let $(P_0,...,P_m)$ be a Clifford system on $\mathbb{R}^{2l}$. We define $m_1:=m,~m_2:=l-m-1$ and
$F:~\mathbb{R}^{2l}\rightarrow \mathbb{R}$ by
$$F(x):=\langle x,x\rangle^2-2\sum_{i=0}^{m}\langle P_ix,x\rangle^2.$$
Then $F$ satisfies the Cartan-M\"{u}nzner differential equations $(2.6)$. If $m_2>0$, then the levels of $F$ form an isoparametric
family with $g=4$ distinct principal curvatures with multiplicities $(m_1, m_2)$. ( Note $m_3=m_1$, $m_4=m_2$ through $2.4$.)
}

\noindent
\emph{Proof}. We have
$$\texttt{grad}_xF=4\langle x,x\rangle x-8\sum\langle P_ix,x\rangle P_ix$$
and with $3.7~(v)$, we have the first equation of $2.6$. Further
\begin{eqnarray*}
\Delta_x~F&=&4(2l+2)\langle x,x\rangle -2\sum\Big(2\langle P_ix,x\rangle\Delta\langle P_ix,x\rangle+2\|\texttt{grad}\langle P_ix,x\rangle\|^2\Big)\\
&=&8(l+1)\langle x,x\rangle-16\langle x,x\rangle (m+1)\\
&=&8(m_2-m_1)\langle x,x\rangle,
\end{eqnarray*}
and so the second equation of $2.6$ is satisfied. The rest of the assertion follows from \cite{Mun}, likewise see $2.5$, $2.6$.
\vspace{2mm}

\noindent
\textbf{4.2~The focal manifolds.} The focal manifolds of an isoparametric family contain very concentrated geometric information.
We study the focal manifolds of our examples in the next theorem and in the following paragraphs of this section. A part of the results
which we obtain is already known from the general theory of isoparametric hypersurfaces.

\noindent
\textbf{Theorem:}\,\,
{\itshape
With the notation from $4.1$, let $f=F|_{S^{2l-1}}$ and $\Sigma:=\Sigma(P_0,...,P_m)$.

\noindent
$(i)$~For $M_-:=f^{-1}(\{-1\})$, we have
$$M_-=\{x\in S^{2l-1}~|~there~ exists~ P\in\Sigma ~with ~x\in E_+(P)\}.$$
In the case $m_2<0$, then $f=-1$, thus $M_-=S^{2l-1}$; this is only possible for $m\in\{1,2,4,8\}$.

In the case $m_2\geq 0$, then $M_-$ is diffeomorphic to the total space of an $(l-1)$-sphere bundle
$$\Gamma:=\{(x,P)~|~x\in S^{2l-1}, ~P\in \Sigma,~x\in E_+(P)\}\xrightarrow{\pi}\Sigma,~~(x,P)\mapsto P.$$
The diffeomorphism from $\Gamma$ onto $M_-$ is furnished by $(x,P)\mapsto x.$ In particular, if $f$
is not constant, then $M_-$ is a -- trivially connected -- submanifold of codimension $m_2+1$ in the
sphere $S^{2l-1}$.

In the case $m_2=0$, then $M_-$ is a hypersurface; this is only possible for $m\in\{1,3,7\}$.

In the case $m_2>0$, then $M_-$ is the focal manifold corresponding to the principal curvatures
of the family of multiplicity $m_2$. The hypersurfaces are $m_2$-sphere bundles over the connected
sphere bundle space $M_-$.

Suppose $(P_0,...,P_m)$ can be extended to a Clifford system $(P_0,...,P_{m+1})$, which by $3.5$ is
not rare, then $\pi:~\Gamma \rightarrow \Sigma$ is trivial and $M_-$ is diffeomorphic to $S^{l-1}\times S^m$.
For $m\equiv 0 (mod~4)$, the geometrically inequivalent Clifford systems (see $3.5$) lead to inequivalent
sphere bundles $\Gamma \rightarrow \Sigma$.
\vspace{1mm}

\noindent
$(ii)$~For $M_+:=f^{-1}(\{+1\})$, we have
$$M_+=\{x\in S^{2l-1}~|~\langle P_0x,x\rangle=...=\langle P_mx,x\rangle=0\}.$$
In the case $m_2\geq 0$, then $M_+$ is a non-empty submanifold of codimension $m_1+1$
and a focal manifold of the level hypersurfaces corresponding to the principal curvatures
of multiplicity $m_1$. The normal bundle of $M_+$ is trivial with $x\rightarrow (P_0x,...,P_mx)$ as
a basis field. Hence the hypersurfaces are trivial sphere bundles over $M_+$.
\vspace{1mm}

\noindent
$(iii)$~For $x\in M_+$ and $P\in \Sigma$, on the normal great circle $c(t):=\cos t~x+\sin t~Px$, we have
$f(c(t))=\cos 4t$.

The normal great circle meets $M_+$ again after $\frac{\pi}{2}$, so that the hypersurface
$M_t=f^{-1}(\{t\})$ of the family in addition to $\cot t$ has also $\cot (t+\frac{\pi}{2})$ as a
principal curvature of multiplicity $m_1$. In the case $m_2=0$, then only two of these principal
curvatures arise: one obtains an
isoparametric family with $g=2$. If $m_2>0$, then the normal great circle meets $M_-$ at
$t=\frac{\pi}{4}$ and $\frac{3\pi}{4}$, so that $\cot(t+\frac{\pi}{4})$ and $\cot(t+\frac{3\pi}{4})$
are principal curvatures with the same multiplicity $m_2$.
}
\vspace{2mm}

\noindent
\emph{Proof}. $(i)$~From $f(x)=-1$, it follows from the definition of $F$ that
$$\langle \Sigma\langle P_ix,x\rangle P_ix, x\rangle=1$$
thus $$\Sigma\langle P_ix,x\rangle P_ix= x,$$
and $P=\Sigma\langle P_ix,x\rangle P_i\in\Sigma$ has $x$ as $+1$-eigenvector. Conversely, if $Px=x$ for an $x\in S^{2l-1}$
and $P\in \Sigma$, then one can assume as in $3.7(iii)$ that $P=P_0$. But then one gets from $3.2(i)$ that
$\langle P_ix,x\rangle=0$ for all $i>0$. It follows that
$$f(x)=1-2\langle P_0x,x\rangle^2=-1$$
and one has the stated characterization of $M_-$.

In the case $m_2<0$, one has $l\leq m$. For $x\in S^{2l-1}$, let $x=x_++x_-$ be the decomposition of $x$ into
eigenvectors of $P_0$, where $x_{\pm}\in E_{\pm}(P_0)$. Then $P_1x_{\pm},...,P_mx_{\pm}$ are orthogonal in $E_{\mp}(P_0)$.
Thus $l=m$, and a simple direct calculation shows that $F(x)=-1$.

Now let $m_2\geq 0$. First of all,
$\pi:~\Gamma \rightarrow\Sigma$ is actually a sphere bundle, since for $P\in\Sigma$, $E_+(P)\times \Sigma \rightarrow \widetilde{\Gamma}$,~
$(x,Q)\mapsto (Id+Q)x$ is a local trivialization of the corresponding vector bundle
$\widetilde{\Gamma}\rightarrow \Sigma$ in a neighborhood of $P$. As we
have already proven that the map $\Gamma\rightarrow M_-$, $(x,P)\mapsto x$ is surjective, it is easily shown that
the map is a submersion. From $3.7(v)$, it is also known that it is injective, and hence it is a diffeomorphism. The
restriction on the values of $m$ is known by $3.5$. Finally, suppose it is possible to extend $(P_0,...,P_m)$, then
$\pi:~\Gamma \rightarrow\Sigma$ is the restriction of the analogous bundle over $\Sigma(P_0,...,P_{m+1})$ to the
``equatorial sphere" $\Sigma$, and is thus trivial. The vector bundle $\widetilde{\Gamma}\rightarrow \Sigma$ has as
a subbundle of $\mathbb{R}^{2l}\times \Sigma\rightarrow\Sigma$ a canonical metric and covariant derivative. The
curvature tensor $\widetilde{R}$ has the form
$$\widetilde{R}(Q_1,Q_2)=\frac{1}{2}Q_1Q_2|_{E_+(P)},$$
for $P\in\Sigma$ and $Q_1,Q_2\in T_P(\Sigma)=(\mathbb{R}P)^{\bot}$.
If $m\equiv 0(mod~4)$, then one has a characteristic $m$-form $\chi$ on $\Sigma$ by the following definition, for
$P\in\Sigma$ and $Q_1,...,Q_m$ an orthonormal basis of $T_P(\Sigma)$:
\begin{eqnarray*}
\chi&:=&\frac{2^{\frac{m}{2}}}{m!}~\texttt{Trace} \Big( \sum_{\sigma}sign(\sigma)\widetilde{R}(Q_{\sigma(1)},Q_{\sigma(2)})\circ\cdots\circ\widetilde{R}(Q_{\sigma(m-1)},Q_{\sigma(m)})\Big)\ast 1\\
&=&\texttt{Trace}(Q_1\cdots Q_m|_{E_+(P)})\ast1\\
&=&\texttt{Trace}(PQ_1\cdots Q_m|_{E_+(P)})\ast1\\
&=&\frac{1}{2}\texttt{Trace}(PQ_1\cdots Q_m)\ast1\\
&=&\frac{1}{2}\texttt{Trace}(P_0\cdots P_m)\ast1.
\end{eqnarray*}

But by $(3.5)$ the (absolute) trace distinguishes between geometrically equivalence classes.
\vspace{1mm}

\noindent
$(ii)$~From the homogeneity of $F$ and $(2.6)$ it follows that $$\|\texttt{grad}~f\|^2=16(1-f^2).$$
For $m_2\geq 0$, then by $(i)$ $f$ is not constant and must assume the value $1$ as maximum:
so $M_+\neq\emptyset$. The remaining statements are immediately clear.
\vspace{1mm}

\noindent
$(iii)$~From $3.7(iii)$, it suffices to consider the case $P=P_0$. Then
\begin{eqnarray*}
F(c(t))&=&1-2\langle \cos t~P_0x+\sin t~x,\cos t~x+\sin t~P_0x\rangle^2\\
& &-2\sum_{i=1}^{m}\langle \cos t~P_ix+\sin t~P_iP_0x,\cos t~x+\sin t~P_0x\rangle^2\\
&=&1-2(2\cos t~\sin t)^2\\
&=& \cos(4t).
\end{eqnarray*}
The remaining statements are again clear.
\vspace{2mm}

\noindent
\textbf{4.3~Table of small multiplicities.} From $3.7(ii)$ and $(iii)$, geometrically equivalent Clifford
systems give congruent isoparametric families. In paragraph $4.6$ we will show a result in the converse direction. In
both cases the small multiplicities play a special role. We thus give a list of our multiplicities $(m_1,m_2)$ from
Theorem $4.1$ as in $3.5$. So $l=k\delta(m)$, $m_1=m$ and $m_2=l-m-1$.

\begin{figure}[h]
\label{table}
\begin{center}
\includegraphics[width=155mm]{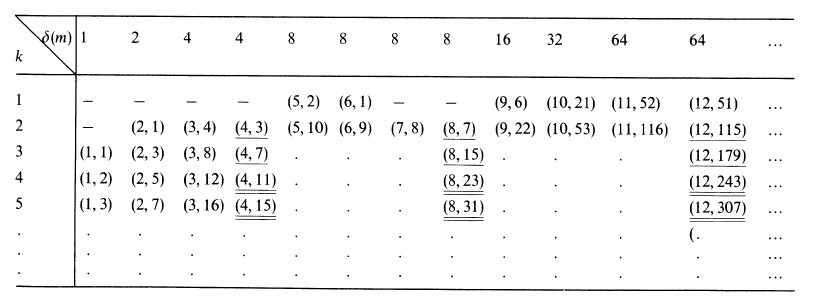}
\end{center}
\end{figure}

The underlined multiplicities \underline{$(m_1,m_2)$} and \underline{\underline{$(m_1,m_2)$}} denote, respectively, the two, \emph{resp.}
three geometrically inequivalent Clifford systems for the multiplicities $(m_1,m_2)$, see $3.5$. We will show further below that all of these with $m\equiv 0 ~mod ~4$ and $l=k\delta(m)$ actually lead to incongruent isoparametric families ( of which there are $[\frac{k}{2}]+1$. )
We will also see that the families for multiplicities $(2,1),~(6,1),~(5,2)$ and one of the two $(4,3)$-families
are congruent to those with multiplicities $(1,2),~(1,6),~(2,5)$ and $(3,4)$, resp. These are all the coincidences under congruence.
\vspace{2mm}

\noindent
\textbf{4.4~Multiplicities of previously known examples.} We give now for comparison the multiplicities of the previously known
examples with $g=4$.
\vspace{1mm}

\noindent
Homogeneous examples:~see \cite{TT}:~$(1,k),~(2,2k-1),~(4,4k-1),~(9,6),~(4,5),~(2,2)$

\noindent
Inhomogeneous examples:~see \cite{OT}:~$(3,4k),~(7,8k).$

We will show that these are all Clifford examples except for the homogeneous examples with $(4,5)$ and $(2,2)$.
\vspace{2mm}

\noindent
\textbf{4.5~The second fundamental tensors of the focal manifolds.}
We now continue the study of focal manifolds which was begun in $4.2$ and will describe their second fundamental tensors.

\noindent
$(i)$~From $4.2(iii)$, we have that for a unit normal $N$ to an $(m_i+1)$-codimensional focal manifold, the corresponding
second fundamental tensor $S_N$ has the eigenvalue $0$ of multiplicity $m_i$ and $+1$, respectively, $-1$ with multiplicities
$m_j$, $\{i,j\}=\{1,2\}$. Thus corresponding to the unit $m_i$-sphere in the normal space, there is an $m_i$-sphere of
symmetric endomorphisms $S$ with $S^3=S$ and $m_i$-dimensional kernel. The classification of these algebraic structures
(these exist likewise for the non-Clifford families with $g=4$) looks much harder than the Clifford case $( S^2=Id )$.
\vspace{1mm}

\noindent
$(ii)$~Let $x\in M_+$. Then by $4.2(ii)$,
$$\perp_xM_+=\{Px~|~P\in\mathbb{R}\Sigma(P_0,...,P_m)\},$$
and for the second fundamental tensors for $N=Px$, $P\in\Sigma:=\Sigma(P_0,...,P_m)$, we have
$$S_Nv=M_+-tangential~component~of~(-Pv).$$
From this and with
$$\Sigma_P:=\{Q\in \Sigma~|~\langle P,Q\rangle=0\},$$
it follows very easily that
\begin{eqnarray*}
&&\ker S_N=\mathbb{R}\Sigma_PN\quad (N=Px)\\
&&S_N|_{(\ker S_N)^{\perp}}=-P|_{(\ker S_N)^{\perp}}.
\end{eqnarray*}
\vspace{1mm}

\noindent
$(iii)$~Let $y\in M_-$ and $P\in\Sigma$ with $Py=y$. Let $\Sigma_P$ again be the equatorial sphere of $\Sigma$
orthogonal to $P$. Since $E_+(P)\cap S^{2l-1}\subset M_-$, from $4.2(i)$, then $\perp_yM_-\subset E_-(P)$.
Further, from $3.7(iii)$, for all $Q\in\Sigma_P$, $F((\cos t~P+\sin t~Q)y)=F(y)=-1$, and hence
$(\cos t~Py+\sin t~Qy)^{'}(0)=Qy\in T_yM_-$. But since $\mathbb{R}\Sigma_Py\subset E_-(P)$ is an $m_1$-dimensional
subspace, it follows for dimensional reasons that
$$\perp_yM_-=\{N\in E_-(P)~|~N\perp\Sigma_Py\}.$$
Now let $N\in\perp_yM_-$, $\|N\|=1$. The normal great circle $\cos t~y+\sin t~N$ meets $M_+$ in $x:=\frac{1}{\sqrt{2}}(y+N)$
and $Px=\frac{1}{\sqrt{2}}(y-N)$.
Since the eigenspaces of the second fundamental tensors of a parallel family of hypersurfaces are parallel along a normal circle,
it follows from $(ii)$ that
$$E_+(S_N)=\mathbb{R}\Sigma_P(y+N),~E_-(S_N)=\mathbb{R}\Sigma_P(y-N)$$
$$\ker S_N=\{v\in E_+(P)~|~v\perp y, ~v\perp \Sigma_PN\}.$$
\vspace{2mm}

\noindent
\textbf{4.6~The uniqueness theorem.} The second fundamental tensors of the focal manifolds are not only well
suited to distinguish families with the same multiplicities as in the following theorems, but also for distinguishing
the geometry at different points of the same focal submanifold, as will become clear further on in our work.

\noindent
\textbf{Theorem:}\,\,
{\itshape
Let $(P_0,...,P_m)$ be a Clifford system on $\mathbb{R}^{2l}$ and
$$m=m_1\leq m_2=l-m-1.$$
Let $f$ be defined as in $4.2$ and $M_-=f^{-1}(\{-1\})$.
Let $P\in\Sigma:=\Sigma(P_0,...,P_m)$ and $y\in E_+(P)\cap S^{2l-1}$.
Then from $4.2(i)$ $y\in M_-$ and $E_+(P)=\mathbb{R}y\oplus Span\cup \ker S_N$, where
the union is taken over all $N\in\perp_yM_-\backslash\{0\}$. Since the Clifford sphere $\Sigma$
is uniquely determined by the set $\{E_+(P)~|~P\in\Sigma\}$, so also is the isoparametric family
uniquely determined. In particular, the congruence class of the hypersurface family determines the
geometric equivalence class of the representation.
}

\vspace{2mm}

\noindent
\emph{Proof.} We denote the orthogonal complement in $T_yM_-$ by $(~~)^{\perp}$.
Then it follows from $4.5(iii)$ that
\begin{eqnarray*}
Span~\displaystyle\bigcup_{N\neq 0}\ker S_N&=&(\displaystyle\bigcap_{N\neq 0}E_+(S_N)\oplus E_-(S_N))^{\perp}\\
&=&(\displaystyle\bigcap_{N\neq 0}\mathbb{R}\Sigma_P(y+N)\oplus \mathbb{R}\Sigma_P(y-N))^{\perp}\\
&=&(\displaystyle\bigcap_{N\neq 0}\mathbb{R}\Sigma_Py\oplus \mathbb{R}\Sigma_PN)^{\perp}\\
&=&(\mathbb{R}\Sigma_Py\oplus \displaystyle\bigcap_{N\neq 0}\mathbb{R}\Sigma_PN)^{\perp},
\end{eqnarray*}
since $\Sigma_Py\subset E_-(P)$ and $\Sigma_PN\subset E_+(P)$.
But since $T_yM_-=\mathbb{R}\Sigma_Py\oplus (\mathbb{R}\Sigma_PN\oplus \ker S_N)$ as
an orthogonal direct sum with $\mathbb{R}\Sigma_PN\oplus \ker S_N\subset E_+(P)$,
then $(\mathbb{R}\Sigma_Py)^{\perp}$ is a subspace of $ E_+(P)$ of dimension
$\dim M_--m_1=m_1+m_2$, so that
$$ E_+(P)=\mathbb{R}y\oplus(\mathbb{R}\Sigma_Py)^{\perp}.$$
Hence, it suffices to prove that
$$\displaystyle\bigcap_{N\neq 0}\mathbb{R}\Sigma_PN=\{0\}.$$
Given a $u\neq 0$ in this intersection, then one has for each $N\neq 0$ a $Q\in \mathbb{R}\Sigma_P$ with
$QN=u$.

From $Q^2=\|Q\|^2Id$, it follows that $N=\frac{1}{\|Q\|^2}Qu$, \emph{i.e. }the rank of the map $\mathbb{R}\Sigma_P \rightarrow \mathbb{R}^{2l}$,
~$Q\mapsto Qu$ is at least $\dim\perp_yM_-=m_2+1$, in contradiction to the assumption that $\dim \mathbb{R}\Sigma_P=m_1\leq m_2$.
\vspace{2mm}

\noindent
\textbf{4.7~Consequences.} A glance at the table $4.3$ shows that the assumption $m_1\leq m_2$ is almost always satisfied.
We will consider the (eight) exceptional cases later. The result $4.6$ is naturally of special interest in the cases $m\equiv 0~mod~4$,
in which there exist $[\frac{k}{2}]+1$ geometrically inequivalent Clifford systems on $\mathbb{R}^{2l}$ with $l=k\delta(m)$.

The corresponding families ( except $\underline{(4,3)}$ and $\underline{(8,7)}$ ) are now seen to be incongruent. If one chooses two hypersurfaces from two
incongruent families at the same distance from $M_-$, then they have the same principal curvatures and thus from the Gauss equation
the same curvature tensor: to two points from the two hypersurfaces, there exists a linear isometry of the corresponding tangent
spaces, which transforms one curvature tensor into the other. Nevertheless, by the remark $2.7$, these hypersurfaces are not
intrinsically isometric!
\vspace{6mm}

\noindent
\textbf{5. Inhomogeneity of Clifford families.} We will later show that several of our examples are homogeneous,
see $6.1$ and $6.3$. Most of our examples are indeed inhomogeneous, since their multiplicities are not found in the list $4.4$
of homogeneous multiplicities. In this section, we will give--without use of the homogeneous classification--a direct geometric
proof of the inhomogeneity of most of the families and their focal manifolds.

Let $(P_0,...,P_m)$ be a fixed Clifford system on $\mathbb{R}^{2l}$. Let the notation $F$, $M$, $M_{\pm}$ etc.be as in $\textbf{4}$.
In particular, let $$\Sigma:=\Sigma (P_0,...,P_m).$$
\vspace{2mm}

\noindent
\textbf{5.1~Theorem:}\,\,
{\itshape
Let $N_+:=\{x\in M_+~|~there~exists~an~orthonormal~Q_0,...,Q_3\in\Sigma~with~Q_0\cdots Q_3x=x\}$. Then $N_+$ has the following
geometric description:
$$N_+=\{x\in M_+~|~there~exists~orthonormal~N_0,...,N_3\in\perp_xM_+~with~\dim(\displaystyle\bigcap_{i}\ker S_{N_i})\geq 3\}.$$
}
\noindent
\emph{Proof.} First let $x\in N_+$, $x=Q_0\cdots Q_3x$. From $4.5(i)$, the vectors $N_i:=Q_ix\in\perp_xM_+$ are
orthonormal and $Q_0Q_1x=-Q_2Q_3x$, $Q_0Q_2x=Q_1Q_3x$ and $Q_0Q_3x=-Q_1Q_2x$ are orthonormal vectors in $\displaystyle\bigcap_{i}\ker S_{N_i}$.
Conversely, let $N_0,N_1\in \perp_xM_+$ with at least three-dimensional intersection of their corresponding kernels. From $4.5(i)$
and $3.7(v)$, there exists an orthonormal $Q_0, Q_1\in\Sigma$ with $Q_ix=N_i$. The intersection of the kernels contains a vector
orthonormal to $Q_0Q_1x$ which by $4.5(i)$ and $3.7(v)$ must be of the form $Q_0Q_2x=Q_1Q_3x$ with $Q_2, Q_3\in\Sigma$ orthonormal and
orthogonal to $Q_0,Q_1$. Hence $Q_0\cdots Q_3x=x$.

\vspace{2mm}

\noindent
\textbf{5.2~Theorem:}\,\,
{\itshape
Suppose $9\leq 3m_1<m_2+9$ and for $m_1=4$ suppose that the additional identity $P_0\cdots P_4\neq \pm Id$ holds, then
$\emptyset\neq N_+\neq M_+$. Thus the focal manifold $M_+$ and the whole family are not homogeneously embedded.
}
\vspace{1mm}

\noindent
\emph{Proof.} The endomorphism $P:=P_0\cdots P_3$ is involutive, symmetric commuting with $P_4,...,P_m$
and anti-commuting with $P_0,...,P_3$. Let $S_+(P):=E_+(P)\cap S^{2l-1}$ be the unit sphere in the $+1$ eigenspace.

For $x\in E_+(P)$, we have
$$F(x)=\langle x,x\rangle^2-2\sum_{i=4}^{m}\langle P_ix,x\rangle^2,$$
since $\langle P_0x,x\rangle=...=\langle P_3x,x\rangle=0$. Thus for $m=m_1=3$, we have that $S_+(P)\subset M_+$.
For $m=4$, it follows from the assumption that $P_4$ is indefinite on $E_+(P)$, and so $S_+(P)\cap M_+=(F|_{S_+(P)})^{-1}(\{1\})$
is $(l-2)$-dimensional. For $m>4$, $P_4,...,P_m$ is a Clifford system on $E_+(P)$, whose $(+)$-focal manifold
of dimension $l-m+2$ is just $S_+(P)\cap M_+$.

In all three cases, thus we have
$$N_+\supset S_+(P)\cap M_+ \neq\emptyset~and ~ \dim(S_+(P)\cap M_+)=l-m+2.$$

On the other hand, by $3.7(iv)$, $E_+(Q_0\cdots Q_3)$ depends only on the orientation of $Span(Q_0,...,Q_3)$. The dimension
of the Grassmann manifold of oriented $4$-planes in $\mathbb{R}\Sigma$ is $4(m+1-4)$. Thus $N_+$ has at most dimension
$$4(m-3)+l-m+2=4(m_1-3)+m_2+3=\dim M_++3m_1-m_2-9<\dim M_+.$$
Thus $N_+\neq M_+$ and the theorem is proven.
\vspace{2mm}

\noindent
\textbf{5.3~Remarks.} The above theorem is not applicable to the case $m\leq 2$ or $m=4$ with $P_0\cdots P_4=\pm Id$,
as these will be shown to be homogeneous. Moreover, because of the assumption $3m_1<m_2+9$, this theorem gives no information
about the further (finite) number of exceptions with $m_1\geq 5$. Thus we wish to generalize the above method in such a way
as to replace $Q_0\cdots Q_3$ in the definition of $N_+$ by products of the form $Q_0\cdots Q_{4\mu-1}$ and $Q_0\cdots Q_{4\mu}$
and their eigenspaces. The geometric interpretation of the eigenvectors will then be more complicated than in Lemma $5.1$.
Occasionally, the geometric interpretation becomes clear in the process of giving the proof, but we will not do the interpretation
in general. Instead we remark here that in the case $m_1\leq m_2$, from $4.6$ and $3.7(iv)$, the Clifford sphere $\Sigma$ and also the set
$E_+(Q_0\cdots Q_{\nu})$ for orthonormal $Q_0,...,Q_{\nu}$ is geometrically determined. Whenever the family of such eigenspaces
meets a level of an isoparametric function in a non-empty proper subset, then that level and hence the whole family
is inhomogeneous.
\vspace{2mm}

\noindent
\textbf{5.4~Inhomogeneous families.} The following theorem is sharper than Theorem $5.2$ when $m_1\geq 5$.

\noindent
\textbf{Theorem:}\,\,
{\itshape
Let $5\leq m_1\leq m_2$. Then the isoparametric family is embedded inhomogeneously.}
\vspace{1mm}

\noindent
\emph{Proof.} The product $P:=P_0\cdots P_4$ is symmetric and involutive; it commutes with $P_0,...,P_4$
and anti-commutes with $P_5,...,P_m$. Thus $(P_0,...,P_4)$ is a Clifford system on the $l$-dimensional eigenspace $E_+(P)$,
and for each hypersurface $M$ of the original family, we have $\dim(M\cap E_+(P))=l-2$. The Grassmann manifold of $5$-planes
in $\mathbb{R}\Sigma$ has dimension $5(m-4)$ and by inspection of $3.5$, we see that $5(m-4)+l-2<2l-2$ since $5\leq m_1\leq m_2$.
Thus $M\cap(\cup E_+(Q_0\cdots Q_4))$, where the union is over all orthonormal $Q_0,...,Q_4$ in $\Sigma$, is a non-empty, proper,
isometry-invariant subset of $M$, and hence $M$ is not homogeneous.
\vspace{2mm}

\noindent
\textbf{5.5~Remarks.} The only Clifford examples with $m_1\geq 5$ for which $m_2<m_1$, and thus to which $5.4$ does
not apply, have multiplicities $(5,2),~(6,1),~(9,6)~or~\underline{(8,7)}$. The first three are homogeneous, and both $\underline{(8,7)}$-families are not (see Section $6$).
\vspace{2mm}

\noindent
\textbf{5.6~Inhomogeneity of the focal manifolds.} If an isoparametric family is inhomogeneous, it is not necessary that its
focal manifolds be inhomogeneous, see example $6.4$ and $6.7(ii)$. But again ``most" focal manifolds are inhomogeneous:

\noindent
\textbf{Theorem:}\,\,
{\itshape
Let $m_1\leq m_2$.

\noindent
$(i)$~If $25\leq 5m_1<m_2+25$, then $M_+$ is inhomogeneously embedded;

\noindent
$(ii)$~If $35\leq 7m_1<m_2+33$, then $M_-$ is inhomogeneously embedded;

\noindent
$(iii)$~$M_+$ in the $(10,21)$-family and $M_-$ in the $(9,22)$-family are inhomogeneously embedded;

\noindent
$(iv)$~If $m_1\equiv 0~mod~4$ and $P_0\cdots P_m\neq \pm Id$, then both focal manifolds are inhomogeneously embedded.}
\vspace{1mm}

\noindent
\emph{Proof.} $(i)$ and $(ii)$. On $E_+(P_0\cdots P_4)$ we have
$$F(x)=\langle x,x\rangle^2-2\sum_{i=0}^{4}\langle P_ix,x\rangle^2,$$
since, in contrast with the proof of Theorem $5.2$, $P_0\cdots P_4$ commutes with $P_0,...,P_4$,
and anti-commutes with $P_5,...,P_m$. The $+1$ eigenspace of the orthogonal five product $Q_0\cdots Q_4$,
for $m_1\geq 5$ cuts out a non-trivial isoparametric family with $m'=4$, $l'=\frac{l}{2}$.
In particular, its intersection with $M_+$, \emph{resp.} $M_-$ is of dimension $4+2(\frac{l}{2}-4-1)=m_1+m_2-5$,\emph{ resp. }
$8+(\frac{l}{2}-4-1)=\frac{1}{2}(m_1+m_2+7)$. In each case, one adds the dimension of the Grassmann manifold of
oriented $5$-planes in $\mathbb{R}\Sigma$, that is $5(m_1-4)$, and compares with the dimension of $M_+$, \emph{resp.} $M_-$,
and the results follow as in $5.2$.

$(iii)$~This follows analogously using the $9$-product $Q_0\cdots Q_8$ and for

$(iv)$~one uses the highest dimensional eigenspace of $P_0\cdots P_m$ in the corresponding way.
\vspace{2mm}

\noindent
\textbf{5.7~Remark on the exceptional cases.} In Section $6$, we will show the homogeneity
of the families with $(m_1,m_2)=(9,6)$, $min\{m_1,m_2\}\leq 2$, and families with $m_1=4$, $P_0\cdots P_4=\pm Id$.
If one takes these results as true, the Table $4.3$ and $5.2$ and $5.6$ show that only for the following cases is the homogeneity
or inhomogeneity of the focal manifolds yet undetermined:

$(i)$~$M_+$ for the multiplicities $\underline{(8,7)}$ and $(8,15)$, the latter with $P_0\cdots P_8=\pm Id$.

$(ii)$~$M_-$ for the multiplicities $(3,4k), (6,9), (7,8), (7,16), \underline{(8,7)}, (10,21), (12, 51)$ and $(8,15)$ and $(8,23)$ in the
case $P_0\cdots P_8=\pm Id$.

In several of these cases, the question of homogeneity remains open, but the question of homogeneity for the whole family for $\underline{(8,7)}$
will be clarified in the following.
\vspace{2mm}

\noindent
\textbf{5.8~The Condition (A) of Takeuchi-Ozeki.} We return now to the beginning of this section. Then (and implicitly likewise
in both of the other proofs in $\textbf{5}$) it was necessary to consider the behavior of the kernel of the second fundamental tensor
for different normal directions $N$ to establish ``inhomogeneous properties". A particularly extreme case occurs, and one which
is of special interest in light of Remark $4.5(i)$, when at some point of a focal manifold the kernels of all second fundamental tensors
$S_N (N\neq 0)$ coincide. We point out here that this is just the condition (A) of Takeuchi and Ozeki which was introduced in Lemma $17$ \cite{OT}.

Concerning Condition (A) in the Clifford case, we have:

\noindent
\textbf{Theorem:}\,\,
{\itshape
Let $m_1\equiv 3 ~mod~4$. For $x\in M_+$, let $d(x):=\dim \cap \ker S_N$ where the intersection is taken over all $N\in \perp_xM_+$. Then we have:
$$\min d= 0,\quad \max d>0.$$
For $m_1\in\{3,7\}$, $\max d=m_1$, \emph{i.e.} there exists a point at which condition (A) is fulfilled.
Condition (A) and the weaker condition $d(x)>0$ are thus ``inhomogeneous properties" on $M_+$.}
\vspace{1mm}

\noindent
\emph{\textbf{Remark.}} Condition (A) holds at certain points of $M_-$ in both $(8,7)$ families (see Section $6.6$).
\vspace{1mm}

\noindent
\emph{Proof.} We first show the second assertion is true. For $m_1=3$, $\max d=3$ by $5.1$ and $5.2$.

For $m_1=7$, choose $x\in S^{2l-1}$ as a common eigenvector of the commuting $4$-products
$$P_0P_1P_2P_3,~P_0P_1P_4P_5,~P_0P_1P_6P_7,~P_0P_2P_4P_6.$$
Since each $P_i$ anti-commutes with at least one of these operators, we get $x\in M_+$. And since $x$ is also an eigenvector of the
product of the above-mentioned operators, one obtains the following identities.

\begin{figure}[h]
\label{identities}
\begin{center}
\includegraphics[width=80mm]{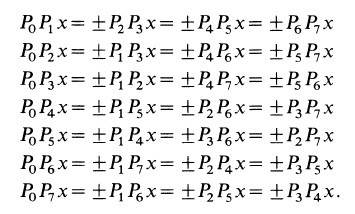}
\end{center}
\end{figure}

From $4.5(ii)$ each of these vectors lies in the kernel of each of the second fundamental tensors $S_{P_ix}$;
from $3.7(v)$, these vectors are also orthonormal. Thus $d(x)=7$.

Now let $m_1=4\mu-1$ for any $\mu$. Choose $x\in S^{2l-1}$ to be a common eigenvector of the commuting operators
$$P_0\cdots P_{m_1}~~~and~~~P_{2i}P_{2i+1}P_{2j}P_{2j+1}, \quad 0\leq i<j \leq 2\mu-1.$$
Then such an $x\in M_+$ and $P_0P_1x=\pm P_{2i}P_{2i+1}x$ lies in the intersection of all the kernels. Thus $max~d>0$.

Finally, if $d(x)>0$, then there exists $Q_0,...,Q_m\in \Sigma$ with
$$P_0Q_0x=P_1Q_1x=\cdots=P_mQ_mx,$$
and $\langle P_i, Q_i\rangle =0$. After modifying $P_1,...,P_m$ by an orthogonal transformation, if necessary, one can assume
that $P_1=Q_0$ and $Q_1=-P_0$. Now it follows for pairwise distinct indices $i,j,k$ that $P_iP_jx=P_kQ_kx$, but
$$0=\langle P_iP_kx, P_iP_jx\rangle = \langle P_iP_kx, P_kQ_kx\rangle = -\langle P_ix, Q_kx\rangle = -\langle P_i, Q_k\rangle. $$
Thus $Q_2,...,Q_m$ are orthogonal to $P_0$ and $P_1$, and by the same argument, one can assume as above that $P_3=Q_2$ and
$P_2=-Q_3$. By repetition of of this last procedure, one obtains
$$P_{2i}P_{2i+1}x=\pm P_{2j}P_{2j+1}x~~~for~all~~~i,j\in\{0,...,2\mu-1\}.$$
Thus $P_0\cdots P_{m_1}x=\pm x$. Thus if $d(x)>0$, then $x\in E_{\pm}(P_0\cdots P_{m_1})$. These eigenspaces depend only on
$\Sigma$ and intersect $S^{2l-1}$
in spheres of dimension $l-1=m_1+m_2<\dim M_+$. Thus there exists $x\in M_+$ with $d(x)=0$.
 \vspace{2mm}

\noindent
\textbf{5.9~The condition (B) of Takeuchi-Ozeki.} With the notation from $5.8$, let $x\in M_+$ and $d(x)=m_1$,
\emph{i.e.} condition (A) holds at $x$. One splits $\mathbb{R}^{2l}=\mathbb{R}x\oplus \perp_xM_+\oplus K\oplus B$ where
$K$, resp. $B$, are the common kernel, resp. the common image, of the second fundamental tensors at $x$, and one splits
$y\in \mathbb{R}^{2l}$ correspondingly as $y=\lambda x+N+k+b$, then one can find $F(y)$ on the individual summands
using the degree of homogeneity. For the term that is constant on $\mathbb{R}x$ and linear on $\perp_xM_+$ one obtains:
$$-8 \sum \langle P_iN,k\rangle \langle P_ib,b\rangle.$$
This is just condition (B) which Takeuchi and Ozeki examined, see \cite{OT}. The main result of
their work is the classification of the $(g=4)$-isoparametric families for which conditions (A) and (B)
hold at a point of a focal manifold. Together with $5.8$, we thus obtain:

\noindent
\textbf{Theorem:}\,\,
{\itshape
The Clifford series of multiplicities $(3,4k)$ and $(7,8k)$ are just the inhomogeneous series of Takeuchi and Ozeki \emph{\cite{OT}}.}
\vspace{4mm}

\noindent
\textbf{6.~Clifford-examples with small $m_1$.}
The arguments of the last section assumed often that $m_1$ was not too small. We will now handle a series of questions
about the exceptional cases--which have yet to be answered.
\vspace{2mm}

\noindent
\textbf{6.1~The homogeneous Clifford series.}

\noindent
\textbf{Theorem:}\,\,
{\itshape
The Clifford isoparametric families of multiplicities $(1,k),~(2,2k-1)$ and the families of multiplicity $(4,4k-1)$ with $P_0\cdots P_4=\pm Id$
are homogeneous.}
\vspace{1mm}

\noindent
\emph{Proof.}~ Let $\mathbb{F}\in \{\mathbb{R}, \mathbb{C}, \mathbb{H}\}$ and $m:=\dim_{\mathbb{R}}\mathbb{F}$.
Let $e_1,...,e_{m-1}$ be the canonical imaginary units of $\mathbb{F}$ and $E_j:~\mathbb{F}^n \rightarrow \mathbb{F}^n$ the left multiplication
by $e_j$. Then as in $3.3$, $E_1,...,E_{m-1}$ induce a Clifford system $(P_0,...,P_m)$ on $\mathbb{R}^{2l}=\mathbb{F}^{2n}$,
where $l=n\cdot m,~n\geq 2$. One sets $c_1:=1,c_j:=e_{j-1}$ for $j\in \{2,...,m\}$ and thus the isoparametric function
for $(P_0,...,P_m)$ is given by
$$F(u,v)=(\|u\|^2+\|v\|^2)^2 - 2\Big\{  (\|u\|^2-\|v\|^2)^2 + 4 \sum\langle u,c_iv\rangle^2 \Big\}~~for ~(u,v)\in \mathbb{F}^n\oplus \mathbb{F}^n.$$
$F$ is invariant under the following sets of isometries.
\begin{eqnarray*}
&&I_1:=\{\cos t~P_0+\sin t~P_1~|~t\in \mathbb{R}\},\\
&&I_2:=\{Id\oplus \alpha Id~|~\alpha \in \mathbb{F},~|\alpha|=1\},\\
&&I_3:=\{A\oplus A~|~A\in U(n,\mathbb{F})\}.
\end{eqnarray*}
The invariance under $I_1$ holds for any Clifford system, see $3.7(iii)$, invariance under $I_2$ and $I_3$ is essentially based
on the special form of $\mathbb{F}$. In particular, the invariance under $I_2$ follows from the fact that the $c_iv$ constitute an
orthogonal basis of $\mathbb{F}v$. Let $(u,v)\in S^{2l-1}$. From the invariance under $I_1$, we can assume without loss of generality
that $\|u\|=\|v\|$. The invariance under $I_3$ further allows us to take without loss of generality $u=(1,0,...,0)/\sqrt{2}$
and using the invariance under $I_2$, take $v=(v_1,...,v_n)$ with $v_1\in (0, \infty)$. Finally, by repeated use of $I_3$,
we get that each $(u,v)\in S^{2l-1}$ lies on the same $F$-level as a point of the form
$$(\bar{u},\bar{v})=\Big((1,0,...,0), (\cos t,\sin t,0,...,0)\Big)/\sqrt{2},~t\in [0,\frac{\pi}{2}].$$
But $F(\bar{u},\bar{v})=-\cos 2t,$ and thus the isometries $I_1, I_2, I_3$ operate transitively on the levels of $F$,
and the family is homogeneous. For $\mathbb{F}= \mathbb{R}, \mathbb{C}, \mathbb{H}$ one gets a family with $m_1=1,2,4,$ where in the last case
$P_0\cdots P_4=-Id$.
\vspace{2mm}

\noindent
\textbf{6.2~The canonical Killing fields.}
We will frequently make use of the following well-known fact: If $P,Q\in \Sigma$ are orthonormal, then $PQ$ is skew-symmetric and
$x \mapsto PQx$
is thus a Killing field on $S^{2l-1}$. But from $3.7(iii)$ we have
$$F((\cos t~Q+\sin t~P)Qx)=F(Qx)$$
and thus the Killing field is tangential to the levels of $F$. By the way, it is easy to show that the vector space spanned by the products
$PQ$ with $P,Q\in \Sigma$ orthonormal is a Lie subalgebra of $\mathfrak{so}(2l)$ isomorphic to $\mathfrak{so}(m+1)$, \emph{i.e.} $Spin(m+1)$
operates on the isoparametric family through isometries.
\vspace{2mm}

\noindent
\textbf{6.3~The homogeneity of the $(9,6)$ family.}
From $4.4$ there exists a homogeneous family with $g=4$ and multiplicities $9$ and $6$, whose order depends on the orientation
so that we cannot say whether $m_1$ is $6$ or $9$.
On the other hand, from $4.3$ there exist Clifford families of multiplicities $(6,9)$ and $(9,6)$. The first is not homogeneous by $5.4$. We show that the other is homogeneous.

\noindent
\textbf{Theorem:}\,\,
{\itshape
The Clifford family with multiplicities $(9,6)$ is homogeneous. Thus, in particular, it is not congruent to the $(6,9)$ family.}

\noindent
\emph{Proof.} For the Clifford system $(P_0,...,P_9)$ on $\mathbb{R}^{32}$, we choose $x\in S^{31}$ to be
a common eigenvector of the commuting operators
$$P_{2i}P_{2i+1}P_{2j}P_{2j+1},\quad 0\leq i< j\leq 4.$$
We set $z(t)=\cos t~x+\sin t~P_0x$ and
\begin{eqnarray*}
&& A(t):=Span\{P_0P_1z(t), P_2P_3z(t)\}\\
&& B(t):=Span\{P_1P_{2i}z(t), P_1P_{2i+1}z(t), P_0P_{2i}z(t), P_0P_{2i+1}z(t)~|~1\leq i\leq 4\}\\
&& C(t):=Span\{P_{2i}P_{2j}z(t), P_{2i}P_{2j+1}z(t)~|~1\leq i<j\leq 4\}
\end{eqnarray*}
Then $A(t), B(t), C(t)$ are subspaces of the tangent space to the $F$-level at $z(t)$ which are pairwise orthogonal to one another
and one finds:
\begin{eqnarray*}
&&\dim A(0)=1,~~\dim A(\frac{\pi}{4})=2\\
&&\dim B(0)\geq 8,~~\dim B(\frac{\pi}{4})=16\\
&&\dim C(0)=12,~~\dim C(\frac{\pi}{4})\geq 6.
\end{eqnarray*}

Moreover, one notices that the generating vectors  in $C(0)$ and $B(\frac{\pi}{4})$ are all orthogonal; in the other cases one finds sufficiently
many (orthonormal) vectors of the form $P_rP_sz(t)$ with fixed $r$. Both inequalities are actually equalities since $z(0)\in M_+^{21}$,
$z(\frac{\pi}{4})\in M_-^{24}$.
The set of $t$-values with $\dim(A(t)\oplus B(t)\oplus C(t))<2+16+12=30$ has only isolated points. Thus one finds a family of hypersurfaces
with $30$ linearly independent Killing fields at a point. Thus the whole family is homogeneous.
\vspace{2mm}

\noindent
\textbf{6.4~The homogeneity of $M_-$ for the $(3,4k)$-families.}
The families of multiplicities $(3,4k)$ are inhomogeneous and by $5.2$ their focal manifold $M_+$
is also inhomogeneous. In contrast, we have

\noindent
\textbf{Theorem:}\,\,
{\itshape
The focal manifolds $M_-$ of the families of multiplicities $(3,4k)$ are homogeneous.}
\vspace{1mm}

\noindent
\emph{Proof.} We will use the notation of the proof of $6.1$ with $\mathbb{F}=\mathbb{H}$ and
consider the Clifford system $(P_1,...,P_4)$. The corresponding function is
$$F(u,v)=(\|u\|^2+\|v\|^2)^2-8\sum_{i=1}^4\langle u,c_iv\rangle^2.$$
For $(u,v)\in M_-$, we have $\|u\|=\|v\|$, and so $F$ is invariant under the set of isometries $I_2$ and $I_3$,
and one obtains the conclusion as in $6.1$.
\vspace{2mm}

\noindent
\textbf{6.5~Coincidence of the Clifford examples.}

\noindent
\textbf{Theorem:}\,\,
{\itshape
The Clifford families with multiplicities $(2,1),~(6,1),~(5,2)$ are congruent to those with multiplicities $(1,2),~(1,6),~(2,5)$.
The $(4,3)$ family with with $P_0\cdots P_4\neq \pm Id$ is congruent to the family with multiplicities $(3,4)$, and thus
from $5.2$ it is inhomogeneous and hence by $6.1$, it is not congruent to the other $(4,3)$ family.}
\vspace{1mm}

\noindent
\emph{Proof.}
For the Clifford system $(P_0,...,P_8)$ on $\mathbb{R}^{16}$ we have from $4.2(i)$, case $m_2<0$,
$$2\sum_{i=0}^8\langle P_ix,x\rangle^2=2\langle x,x\rangle^2$$
and so
$$\langle x,x\rangle^2-2\sum_{i=0}^k\langle P_ix,x\rangle^2=-\Big(\langle x,x\rangle^2-2\sum_{i=k+1}^8\langle P_ix,x\rangle^2\Big).$$
For $k\in\{4,5,6\}$, $P_0,...,P_k$ and $P_{k+1},...,P_8$ are Clifford systems on $\mathbb{R}^{16}$ and in the case $k=4$ we have
$P_0\cdots P_4\neq \pm Id$, since this product anti-commutes with $P_5$. Thus the result follows in the cases with $m_1+m_2=7$.
The $(2,1)$-case can be proven analogously.
\vspace{2mm}

\noindent
\textbf{6.6~Families with multiplicities \underline{$(8,7)$}.}
We will call the family with $P_0\cdots P_8=\pm Id$, where without loss of generality we take the $+$ sign, the \emph{definite family};
the other we will call the \emph{indefinite family}.

\noindent
\textbf{Theorem:}\,\,
{\itshape
$(i)$~For the indefinite $(8,7)$-family, both focal manifolds are inhomogeneously embedded.

\noindent
$(ii)$~For the definite $(8,7)$-family, $M_+$ is homogeneously embedded, $M_-$ inhomogeneously embedded. This and $(i)$
imply that the two $(8,7)$ families are not congruent to one another.

\noindent
$(iii)$~Neither $(8,7)$ family is congruent to the $(7,8)$ family.}
\vspace{1mm}

\noindent
\emph{Proof.}
$(i)$~From the classification in $3.5$, we can extend a Clifford system $(P_0,...,P_8)$ on $\mathbb{R}^{32}$ to a system $(P_0,...,P_9)$,
and we can use the considerations from the proof of $6.3$. Let $x\in S^{31}$ again be a common eigenvector of $P_{2i}P_{2i+1}P_{2j}P_{2j+1}$,
$0\leq i<j\leq 4$. Then, in particular, $P_{2i}P_9x=\pm P_{2i+1}P_8x$ and the right side is a Killing field of our $(8,7)$-family.
Thus, with the same definitions as in $6.3$, one obtains $21$ Killing fields on $M_+$, that all lie in the span of the kernels of the second fundamental tensors, see $4.5(ii)$, and at $x$ they span a $21$-dimensional space $A(0)\oplus B(0)\oplus C(0)$. Thus with
$$\sigma(x):=\dim Span\{v\in T_xM_+~|~there~exists~N\in\perp_xM_+\setminus \{0\},~S_Nv=0\}$$
we have $\sigma(x)\geq 21$.

On the other hand, we choose $u\in S^{31}$ to be a $(+1)$-eigenvector of $P_9$, so that by $4.2(ii)$, $u\in M_+$. Further,
we have $P_9(P_iP_ju)=P_iP_ju$ for all $i,j\in\{0,...,8\}$. Thus from $4.5(ii)$
$$\sigma(u)\leq \dim E_+(P_9)=16.$$
Thus $M_+$ is inhomogeneously embedded. Note that for $x$ as chosen above we have $\sigma(x)=21=\dim M_+-1$, and so
the isometry group of $M_+$ has orbits of codimension $1$.

To study $M_-$, we first choose $x$ as above. Then $y=(x+P_0x)/\sqrt{2}\in M_-$ with $P_0y=y$. Obviously $y$
is likewise an eigenvector of the $4$-products $P_{2i}P_{2i+1}P_{2j}P_{2j+1}$ with $0\leq i<j\leq 4$.
With the help of $3.7(iv)$, one verifies for $1<i<j\leq 9$ and $1\leq k\leq 9$ that
\begin{eqnarray*}
&&\langle P_1P_iP_jy,P_ky\rangle=0\\
&&\langle P_0P_9y,P_ky\rangle=0,\quad k\neq 9.
\end{eqnarray*}
Thus, one sees with $4.5(iii)$ that for a fixed $i\in\{2,...,9\}$
$$\perp_yM_-=Span\{P_0P_9y, P_1P_iP_jy~|~2\leq j\leq 9,~j\neq i \},$$
and the $8$ listed vectors make an orthonormal basis. Further, it follows from $4.5(iii)$ that
\begin{eqnarray*}
&& P_1P_iy \in \ker S_{P_0P_9y},~~2\leq i\leq 8\\
&& P_2P_ry \in \ker S_{P_1P_2P_jy},~~r,j\in\{3,...,9\},~~r\neq j,
\end{eqnarray*}
and that these $14$ vectors are pairwise orthogonal.
With $\sigma$ defined as above for $M_-$ instead of $M_+$, we thus have $\sigma(y)\geq 14$.
(By the way, one can show that $P_1P_9y$ is perpendicular to all the kernels, and thus $\sigma(y)=14$. As for $M_+$,
one can further show that on the $23$-dimensional $M_-$, there are $22$ Killing fields that are linearly independent at $y$.)

\noindent
\emph{Translator's Note}: Regarding the statement in the parentheses above, by the last line above $4.6$ one gets that
$$Span\Big(\ker(S_N)\Big)\subset E_+(P_0)\cap y^{\perp},$$
which has dimension $15$. Since the vector $P_1P_9y$ is in $E_+(P_0)\cap y^{\perp}$ and is perpendicular
to all the kernels, one gets $\sigma(y)\leq 14$. Since $\sigma(y)\geq 14$ has already been shown, it follows that $\sigma(y)= 14$.
\vspace{5mm}

Finally, we show that there is a point $v$ on $M_-$ with $\sigma(v)<14$.
Let $v\in S^{31}\cap E_+(P_0)\cap E_+(P_1\cdots P_8)$. Since $P_0$ and $P_1\cdots P_8$ commute,
but $P_9$ and $P_0$ anti-commute while $P_9$ and $P_1\cdots P_8$ commute, there exists such a $v$ in
$M_-$ by $4.2(i)$. Further it follows that the eigenspaces of $P_0$ and $P_1\cdots P_8$ have $8$-dimensional
intersection. For $N\in E_-(P_0)\cap E_+(P_1\cdots P_8)$ we have
$$\langle N, P_kv\rangle=0,\quad 1\leq k\leq 8,$$
since $P_kv\in E_-(P_0)\cap E_-(P_1\cdots P_8)$, as $P_k$ anti-commutes with both operators. Thus on dimensional grounds,
$$\perp_vM_-=E_-(P_0)\cap E_+(P_1\cdots P_8).$$

For $N\in \perp_vM_-\backslash \{0\}$, we have
$$\ker S_N=\{w\in E_+(P_0)~|~0=\langle w,v\rangle=\langle w,P_1N\rangle=...=\langle w,P_8N\rangle\}.$$
But since $P_iN\in E_+(P_0)\cap E_-(P_1\cdots P_8)$, it follows on dimensional grounds that
$$\ker S_N=\{w\in E_+(P_0)\cap E_+(P_1\cdots P_8)~|~\langle w,v\rangle=0\}.$$
Since the kernel of $N$ is not dependent on $N$, condition (A) holds at $v$ and
$\sigma(v)=7.$ Thus $M_-$ is inhomogeneous.

$(ii)$~Without loss of generality, we can take $P_0\cdots P_8=+Id$. We first show that there exists
$$x\in M_+\cap E_+(P_2P_4P_6P_8)\cap E_+(P_3P_4P_7P_8)\cap E_+(P_5P_6P_7P_8).$$

Since $P_2P_4P_6P_8$ anti-commutes with $P_2$, $E_+(P_2P_4P_6P_8)$ has dimension $16$. It is an invariant subspace of the
anti-commuting operators $P_3P_4P_7P_8$ and $P_3$. Thus $E_+(P_2P_4P_6P_8)\cap E_+(P_3P_4P_7P_8)$ is of dimension $8$
and further it is an invariant subspace of the
anti-commuting operators $P_5P_6P_7P_8$ and $P_5$. Thus $E_+(P_2P_4P_6P_8)\cap E_+(P_3P_4P_7P_8)\cap E_+(P_5P_6P_7P_8)$
is of dimension $4$ and on this space, we have
$$F(x)=\langle x,x\rangle^2-2\sum_{i=0}^1\langle P_ix,x\rangle^2.$$

This function is not constant and a maximum point lies in $M_+$. We choose such an $x$. This is then an eigenvector of
$P_3P_4P_5P_6$, and since $P_0\cdots P_8=Id$, of $P_0P_1P_2P_{2i-1}P_{2i}$, $i\in \{2,3,4\}$. Thus as in the proof of $5.8$,
one obtains the following identities.

\begin{figure}[h]
\label{identities1}
\begin{center}
\includegraphics[width=110mm]{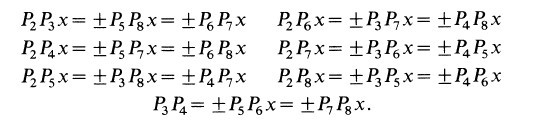}
\end{center}
\end{figure}

These vectors are obviously orthonormal, and with $4.2(ii)$ and $2.7(iv)$ one shows that these can be completed to
form an orthonormal basis of $T_xM_+$ by adjoining
$$P_0P_1x,...,P_0P_8x,P_1P_2x,...,P_1P_8x.$$
Since the vectors are all values at $x$ of Killing vector fields, it follows that $M_+$ is a homogeneous submanifold.

To study $M_-$, we first choose $y\in S^{31}$ as common eigenvector of the commuting $4$-products
$$P_0P_1P_2P_3,~P_0P_1P_4P_5,~P_0P_1P_6P_7,~P_0P_2P_4P_6,$$
as in the proof of $5.8$. Then without restriction $P_8y=y$ and $y\in M_-$.

Let $N\in \perp_yM_-$, $k\in \{0,...,7\}$
and $i\in\{1,...,7\}$. Then from the identities in $5.8$, there exists a $j\in \{1,...,7\}$ with
$$P_0P_iy=\pm P_kP_jy.$$
Thus $\langle P_0P_iy, P_kN\rangle=\pm \langle P_jy, N\rangle=0,$
and so $P_0P_1y,...,P_0P_7y\in \ker S_N$, and the $7$-dimensional kernel is independent of $N\neq 0$, and condition
(A) holds at $y$: so $$\sigma(y)=7.$$

To find a point $v\in M_-$ with other behavior of the kernel, we first establish that from $3.5$ there exist endomorphisms
$\tilde{P_0},\tilde{P_9}$ such that $\{\tilde{P_0}, P_1,...,P_8,\tilde{P_9}\}$ is a Clifford system. We
then show that there exists a point
$$v\in S^{31}\cap E_+(P_0)\displaystyle\bigcap_{1}^4E_+(\tilde{P_0},\tilde{P_9}P_{2i-1}P_{2i}).$$
The $16$-dimensional eigenspace $E_+(P_0)$ is invariant under the anti-commuting operators $\tilde{P_0}$
and $\tilde{P_0}\tilde{P_9}P_1P_2$ ( note $P_0=P_1\cdots P_8$ ).
Thus $E_+(P_0)\cap E_+(\tilde{P_0}\tilde{P_9}P_1P_2)$ is $8$-dimensional and further it is an invariant
subspace of the anti-commuting operators $\tilde{P_0}\tilde{P_9}P_3P_4$ and $P_4P_5$. Thus $E_+(P_0)\cap E_+(\tilde{P_0}\tilde{P_9}P_1P_2)\cap E_+(\tilde{P_0}\tilde{P_9}P_3P_4)$ has dimension $4$ and
is invariant under the anti-commuting operators $\tilde{P_0}\tilde{P_9}P_5P_6$ and $P_6P_7$.
Thus $E_+(P_0)\cap\cdots \cap E_+(\tilde{P_0}\tilde{P_9}P_5P_6)$ has dimension $2$.
But from $P_0=P_1\cdots P_8$, it follows that it is likewise contained in $E_+(\tilde{P_0}\tilde{P_9}P_7P_8)$.
Thus there exists such a $v$, which then naturally lies in $M_-$. For such a $v$, it follows that $\langle \tilde{P_0}P_iv,P_kv\rangle=0$
for all $i,k\in \{1,...,8\}$. Thus from $4.5(iii)$ we have
$$\perp_vM_-=Span\{\tilde{P_0}P_1v,...,\tilde{P_0}P_8v\},$$
and in particular, $P_5\tilde{P_0}P_1v\in (\ker S_{\tilde{P_0}P_1v})^{\perp}$.

We show, on the other hand, that this vector lies in $\ker S_{\tilde{P_0}P_3v}$. Then the kernels are
not all the same at $v$ and $\sigma(v)>7$.

We must from $4.5(iii)$ show that
$$\langle P_5\tilde{P_0}P_1v, P_k\tilde{P_0}P_3v\rangle=\langle P_kP_1P_3P_5v,v\rangle=0,$$
for all $k\in \{1,...,8\}$. But by the choice of $v$, we have
$$\langle P_kP_1P_3P_5v,v\rangle=-\langle P_kP_2\tilde{P_0}\tilde{P_9}P_4\tilde{P_0}\tilde{P_9}P_6\tilde{P_0}\tilde{P_9}v, v\rangle=\langle P_kP_2P_4P_6\tilde{P_0}\tilde{P_9}v, v\rangle=0,$$
from $3.7(iv)$ in the cases $k\not \in \{2,4,6\}$. But if $k=2$, then
$$\langle P_2P_1P_3P_5v, v\rangle= \langle P_3P_5P_2P_1v, v\rangle =\langle P_3P_5\tilde{P_9}\tilde{P_0}v, v\rangle=-\langle P_3P_6v, v\rangle=0,$$
likewise from $3.7(iv)$. One shows this analogously for $k=4,6.$ Thus $\sigma(v)>7$ and $M_-$ is inhomogeneously embedded.

$(iii)$~Were an $(8,7)$-family congruent to a $(7,8)$-family, this would give Clifford systems $(P_0,...,P_8)$ and $(P_9,...,P_{16})$ on
$\mathbb{R}^{32}$ with
$$\langle x,x\rangle^2-2\sum_{0}^8\langle P_ix,x\rangle^2=-\Big(\langle x,x\rangle^2-2\sum_{9}^{16}\langle P_ix,x\rangle^2\Big)$$
or $$\sum_{0}^{16}\langle P_ix,x\rangle^2=\langle x,x\rangle^2.$$

We now show that $(P_0,...,P_{16})$ is a Clifford system on $\mathbb{R}^{32}$, in contradiction to $3.5$. Through differentiation, one obtains
$$\sum_{0}^{16}\langle P_ix,x\rangle P_ix=\langle x,x\rangle x,$$
and for $i\in \{0,...,8\}$ and $u\in E_{\pm}(P_i)$
$$\langle u,u\rangle u+\sum_{9}^{16}\langle P_iu,u\rangle P_iu=\langle u,u\rangle u,$$
thus $$\langle P_9u,u\rangle=...=\langle P_{16}u,u\rangle=0.$$
Thus for $x=u+v$ with $u\in E_+(P_i)$, $v\in E_-(P_i)$ and $j\in \{9,...,16\}$
$$\langle P_iP_jx,x\rangle=\langle P_j(u+v), P_i(u+v)\rangle=\langle P_ju,u\rangle-\langle P_jv,v\rangle=0,$$
and so $P_iP_j+P_jP_i=0$, and thus $(P_0,...,P_{16})$ is a Clifford system. Contradiction!
\vspace{2mm}

\noindent
\textbf{6.7~Summary.} The results $4.6,~6.1,~6.3,~6.5$ and $6.6$
show:
The congruence classes of the families with multiplicities $(1,2),~(1,6),~(2,5),~(3,4)$ each occur a second
time in the Table $4.3$, all others are listed only once.

\noindent
Received on 24. September 1980
\vspace{8mm}

\noindent
\textbf{Added in proof}

\noindent
We just now received:

\noindent
Dorfmeister, J., Neher, E.: {\it Isoparametric triple systems of algebra type}, preprint.
\vspace{2mm}

\noindent
A major result is: In the case $g=4$, it follows from Condition (A) that the isoparametric
families are either Clifford or that $(m_1, m_2)=(2, 2)$.

\end{document}